

\baselineskip=14pt
\parskip=10pt

\font\eightrm=cmr8 

\magnification=\magstephalf
\def\P{{\cal P}}

\def\1{{\overline{1}}}
\def\2{{\overline{2}}}
\parindent=0pt
\overfullrule=0in

\def\frac#1#2{{#1 \over #2}}
\centerline
{\bf An Experimental Mathematics Approach to the Area Statistic of Parking Functions  }
\bigskip
\centerline
{\it Yukun YAO and Doron ZEILBERGER}
\bigskip
{\bf Abstract.} We illustrate the experimental, empirical, approach to  mathematics (that contrary to popular belief, is often rigorous),
by using {\it parking functions} and their  `area' statistic, as a {\it case study}.
Our methods are purely {\it finitistic} and elementary, taking full advantage, of course, of our beloved silicon  servants. 

\bigskip

{\bf Accompanying Maple package and input and output files}

This  article is accompanied by a Maple package {\tt ParkingStatistics.txt} available from the front of this article

{\tt http://sites.math.rutgers.edu/\~{}zeilberg/mamarim/mamarimhtml/par.html} \quad , 

where readers can also find lots of output files, and nice pictures.

{\bf Pre-History (and Pre PC)}

Once upon a time, way back in the nineteen-sixties, there was a one-way street (with no passing allowed), with $n$  parking spaces
bordering the sidewalk. Entering the street were  $n$ cars, each driven by a loyal {\it husband}, and sitting next to him,
dozing off, was his capricious (and a little bossy) {\it wife}. At a random time (while still along the street), the wife wakes
up and orders her husband, {\bf park here, darling!}. If that space is unoccupied, the hubby gladly obliges,
and if the parking space is {\bf occupied}, he  parks, if possible, at the first still-empty parking space.
Alas, if all the latter parking spaces are occupied, he has to go around the block, and drive back to the
beginning of this one-way street, and then look for the first available spot. Due to construction, this
wastes half an hour, making the wife very cranky. 

{\bf Q}: What is the probability that no one has to go around the block?

{\bf A}: $(n+1)^{n-1}/n^n \, \asymp  \, \frac{e}{n+1}$.

Both the question and its elegant answer are due to Alan Konheim and Benji Weiss [KW].

{\bf Parking Functions}

Suppose wife $i$ ($1 \leq i \leq n$) prefers parking-space $p_i$, then the preferences of the wives can
be summarized as an array $(p_1, \dots, p_n)$, where $1 \leq p_i \leq n$. So altogether there are $n^n$ possible
preference-vectors, starting from $(1, \dots , 1)$ where it is clearly possible for everyone  to park,
and ending with $(n,...,n)$ (all $n$), where every wife prefers the last parking space, and of course
it is impossible. Given a preference vector $(p_1, \dots, p_n)$, let $(p_{(1)}, \dots, p_{(n)})$
be its {\it sorted} version, arranged in {\it (weakly) increasing order}. \hfill\break
For example if $(p_1,p_2,p_3,p_4)=(3,1,1,4)$ then  $(p_{(1)},p_{(2)},p_{(3)},p_{(4)})=(1,1,3,4)$.

We invite our readers to convince themselves that a parking-space preference vector  $(p_1, \dots, p_n)$
makes it possible for every husband to park without inconveniencing his wife if and only if $p_{(i)} \leq i$ for
$1 \leq i \leq n$. This naturally leads to the following definition.

{\bf Definition of a Parking Function}: A vector of positive integers $(p_1, \dots, p_n)$ with $1 \leq p_i \leq n$ is
a {\bf parking function} if its (non-decreasing) sorted version $(p_{(1)}, \dots, p_{(n)})$ 
(i.e. $p_{(1)} \leq p_{(2)} \leq \dots  \leq p_{(n)}$, and the latter is a permutation of the former)
satisfies
$$
p_{(i)} \leq i \quad, \quad (1 \leq i \leq n) \quad .
$$

As we have already mentioned above, Alan Konheim and Benji Weiss ([KW]) were the {\bf first} to state and {\bf prove} the following theorem.

{\bf The Parking Function Enumeration Theorem}: There are $(n+1)^{n-1}$ parking functions of length $n$.

There are {\bf many} proofs of this lovely theorem, possibly the slickest is due to the brilliant human
Henry Pollak, (who apparently did not deem it worthy of publication. It is quoted, e.g. in [FR]).
It is nicely described on pp. 4-5 of [St1] (see also [St2]), hence we will {\bf not} repeat it here. 
Instead, as a warm-up to the `statistical' part, and to illustrate the power of experiments, we will
give a much uglier proof, that, however, is {\it motivated}.

Before going on to present {\it our} (very possibly not new) `humble' proof, we should mention that one natural way to prove the
Konheim-Weiss theorem is by a {\it bijection} with labeled trees on $n+1$ vertices, that Arthur
Cayley famously proved is also  enumerated by $(n+1)^{n-1}$. The first such bijection, as far as we know,
was given by the great formal linguist, Marco Sch\"utzenberger ([Sc]). This was followed by an elegant bijection
by the {\it classical} combinatorial giants Dominique Foata and John Riordan [FR], and others.

Since we know (at least!) $16$ different proofs of Cayley's formula (see, e.g. [Z3]), and at least four  different bijections
between parking functions and labeled trees, there
are at least $64$ different proofs (see also [St3], ex. 5.49) of the Parking Enumeration theorem. To these one must add proofs
like Pollak's, and a few other ones.

Curiously, our `new' proof has some resemblance to the very first one in [KW], since they both
use {\it recurrences} (one of the greatest tools in the experimental mathematician's tool kit!),
but our proof is (i) motivated (ii) experimental (yet fully rigorous).

{\bf An Experimental Mathematics Motivated Proof of the Kohnheim-Weiss Parking Enumeration Theorem}

When encountering a new combinatorial family, the first task is to write a computer program to
enumerate as many terms as possible, and hope to {\it conjecture} a nice formula.
One can also try and "cheat" and use the great OEIS, to see whether anyone came up with this sequence
before, and see whether this new combinatorial family is mentioned there.

A very brute force approach, that will not go very far (but would suffice to get the first five terms needed for the OEIS)
is to list the {\it superset}, in this case all the $n^n$ vectors in $\{1 \dots n\}^n$ and for each of them sort it,
and see whether the condition $p_{(i)} \leq i$ holds for all $1 \leq i \leq n$. Then count the vectors that
pass this test.

But a much better way is to use {\bf dynamical programming} to express the desired sequence, let's call it $a(n)$,
in terms of values $a(i)$ for $i<n$. 

Let's analyze the anatomy of a typical parking function of length $n$.
A natural parameter is the number of $1$'s that show up, let's call it $k$ ($0 \leq k \leq n$).
i.e.
$$
p_{(1)}=1 \quad ,\quad \dots \quad , \quad  p_{(k)}=1 \quad , \quad 2 \leq p_{(k+1)} \leq k+1 
\quad , \quad \dots  \quad , \quad  p_{(n)} \leq n  \quad .
$$
Removing the $1$'s yields a shorter weakly-increasing vector
$$
2 \leq p_{(k+1)}  \leq p_{(k+2)} \leq  \dots \quad  \leq \, p_{(n)} \quad ,
$$
satisfying  
$$
p_{(k+1)} \leq k+1 \quad , \quad p_{(k+2)} \leq k+2 \quad , \quad \dots \quad , \quad  p_{(n)} \leq n \quad .
$$
Define
$$
(q_1, \dots, q_{n-k})\, := \, (p_{(k+1)}-1, \dots , p_{(n)}-1 ) \quad .
$$
The vector $(q_1, \dots, q_{n-k})$ satisfies
$$
1 \leq q_1 \leq \dots \leq q_{n-k}  \quad ,
$$
and
$$
q_1 \leq k \quad , \quad q_2 \leq k+1 \quad , \quad \dots \quad , \quad q_{n-k} \leq n-1 \quad .
$$

We see that the set of parking functions with exactly $k$ $1$'s may be obtained by taking the above set of vectors of length $n-k$,
adding $1$ to each component, scrambling it in everywhich way, and inserting the $k$ $1$'s in everywhich way.

Alas, the `scrambling' of the set of such $q$-vectors is not of the original form. We are forced to consider a
more {\bf general} object, namely scramblings  of vectors of the form $p_{(1)} \leq \dots \leq  p_{(n)}$ with
the condition
$$
p_{(1)} \leq a \quad , \quad p_{(2)} \leq a+1 \quad , \quad \dots \quad , \quad p_{(n)} \leq a+n-1 \quad ,
$$
for a {\it general}, positive integer $a$, not just for $a=1$.
So in order to get the dynamical programming recurrence rolling we are {\bf forced} to introduce a more general object,
called an {\bf $a$-parking function}. This leads to the following definition.

{\bf Definition of an a-Parking Function}: A vector of positive integers $(p_1, \dots, p_n)$ with $1 \leq p_i \leq n+a-1$ is
an {\bf $a$-parking function} if its (non-decreasing) sorted version $(p_{(1)}, \dots, p_{(n)})$ 
(i.e. $p_{(1)} \leq p_{(2)} \leq \dots  \leq p_{(n)}$, and the latter is a permutation of the former)
satisfies
$$
p_{(i)} \leq a+i-1 \quad, \quad (1 \leq i \leq n) \quad .
$$

Note that the usual parking functions are the special case $a=1$. So if we would be able to find an efficient
recurrence for counting $a$-parking functions, we would be able to answer our original question.

So let's redo the above `anatomy' for these more general creatures, and hope that the two parameters $n$ and $a$ would
suffice to establish a {\bf recursive scheme}, and we won't need to introduce yet more general creatures.

Let's analyze the anatomy of a typical $a$-parking function of length $n$.
Again, a natural parameter is the number of $1$'s that show up, let's call it $k$ ($0 \leq k \leq n$).
i.e.
$$
p_{(1)}=1 \quad ,\quad \dots \quad , \quad  p_{(k)}=1 \quad , \quad  2 \leq p_{(k+1)} \leq a+k \quad , \quad \dots  \quad  p_{(n)} \leq a+n-1  \quad .
$$
Removing the $1$-s yields a sorted vector
$$
2 \leq p_{(k+1)}  \leq p_{(k+2)} \leq  \dots \, \leq \, p_{(n)} \quad ,
$$
satisfying  
$$
p_{(k+1)} \leq k+a \quad , \quad p_{(k+2)} \leq k+a+1 \quad , \quad  \dots \quad , \quad  p_{(n)} \leq n+a-1 \quad .
$$
Define
$$
(q_1, \dots, q_{n-k})\, := \, (p_{(k+1)}-1 \quad , \quad \dots  \quad , \quad p_{(n)}-1 ) \quad .
$$
The vector $(q_1, \dots, q_{n-k})$ satisfies
$$
q_1 \leq \dots \leq q_{n-k}  \quad
$$
and
$$
q_1 \leq k+a-1 \quad , \quad q_2 \leq k+a \quad , \quad \dots \quad , \quad q_{n-k} \leq n+a-1 \quad .
$$

We see that the set of $a$-parking functions with exactly $k$ $1$'s may be obtained by taking the above set of vectors of length $n-k$,
adding $1$ to each component, scrambling it in everywhich way, and inserting the $k$ $1$'s in everywhich way.

But {\bf now} the set of scramblings of the vectors $(q_1, \dots q_{n-k})$ is an {\bf old friend!}. It is the set of $(a+k-1)$-parking functions of
length $n-k$. To get all $a$-parking functions of length $n$ with exactly $k$ ones we need to take each and every
member of the set of $(a+k-1)$-parking functions of length $n-k$, add $1$ to each component, and insert $k$  ones in every which
way. There are ${{n} \choose {k}}$ ways of doing it. Hence the number of $a$-parking functions of length $n$ with exactly $k$ ones
is  ${{n} \choose {k}}$ times the number of $(a+k-1)$-parking functions of length $n-k$.
Summing over all $k$ between $0$ and $n$ we get the following recurrence.

{\bf Fundamental Recurrence for $a$-parking functions}

Let $p(n,a)$ be the number of $a$-parking functions of length $n$. We have the recurrence
$$
p(n,a) \, = \, \sum_{k=0}^{n} \, {{n} \choose {k}} p(n-k,a+k-1) \quad ,
\eqno(FundamentalRecurrence)
$$

subject to the  {\bf boundary conditions} $p(n,0)=0$ for $n \geq 1$, and $p(0,a)=1$ for $a \geq 0$.

Note that in the sense of Wilf [W], this already {\bf answers} the enumeration problem to compute $p(n,a)$ and hence
$p(n,1)=p(n)$, since this gives us a {\it polynomial time} algorithm to compute $p(n)$ (and $p(n,a)$).

Moving the term $k=0$ from the right to the left, and denoting $p(n,a)$ by $p_n(a)$ we have
$$
p_n(a)- p_n(a-1) \, = \, \sum_{k=1}^{n} \, {{n} \choose {k}} p_{n-k}(a+k-1) \quad .
$$

Hence we can express $p_n(a)$ as follows, in terms of $p_{m}(a)$ with $m<n$.
$$
p_n(a) = \sum_{b=0}^{a} \left ( \sum_{k=1}^{n} \, {{n} \choose {k}} p_{n-k}(b+k-1) \right ) \quad .
$$

Here is the  Maple code that implements it

{\obeylines
{\tt
p:=proc(n,a) local k,b:
if n=0 then  
RETURN(1) 
else 
factor(subs(b=a,sum(expand(add(binomial(n,k)*subs(a=a+k-1,p(n-k,a)),k=1..n)),a=1..b))): 
fi:
end:
}
}

If you copy-and-paste this onto a Maple session, as well as the line below,

{\tt [seq(p(i,a),i=1..10)];}

you would {\tt immediately} get

$$
[a,a \left( a+2 \right) ,a \left( a+3 \right) ^{2},a \left( a+4 \right) ^{3},a \left( a+5 \right) ^{4},a \left( a+6 \right) ^{5},a \left( a+7 \right) ^{6},a
 \left( a+8 \right) ^{7},a \left( a+9 \right) ^{8},a \left( a+10 \right) ^{9}] \quad .
$$

Note that these are {\it rigorously} proved exact expressions, in terms of {\it general} $a$ (i.e. {\it symbolic} $a$) for $p_n(a)$, for
$1 \leq n \leq 10$, and we can easily get more. The following {\bf guess} immediately comes to mind

$$
p(n,a)=p_n(a)= a(a+n)^{n-1} \quad .    
$$

How to prove this rigorously? If you set $q(n,a):=a(a+n)^{n-1}$, since $q(n,0)=0$ and $q(0,a)=1$, the fact that
$p(n,a)=q(n,a)$ would follow {\bf by induction} once you prove that $q(n,a)$ also satisfies the same fundamental recurrence.

$$
q(n,a) \, = \, \sum_{k=0}^{n} \, {{n} \choose {k}} q(n-k,a+k-1) \quad .
\eqno(FundamentalRecurrence') 
$$
In other words, in order to prove that $p(n,a)=a(n+a)^{n-1}$, we have to prove the {\bf identity}
$$
a(a+n)^{n-1} \, = \, \sum_{k=0}^{n} \, {{n} \choose {k}} (a+k-1)(a+n-1)^{n-k-1} \quad ,
$$
but this is an {\bf immediate} consequence of the {\bf binomial theorem}, hence trivial to both humans and machines.

We have just rigorously reproved, via experimental mathematics, the following well-known theorem.

{\bf Theorem}: The number of $a$-parking functions of length $n$ is
$$
p(n,a)=a\,(a+n)^{n-1} \quad .
$$
In particular, by substituting $a=1$, we reproved the original Konheim-Weiss theorem that $p(n,1)=(n+1)^{n-1}$.

{\bf From Enumeration to Statistics in General}

Often in enumerative combinatorics, the class of interest has natural `statistics', like height, weight, and IQ for humans, and
one is interested rather than, for a finite set $A$,
$$
|A| \, := \, \sum_{a \in A} 1 \quad,
$$
called the {\it naive counting}, and getting a {\bf number} (obviously a non-negative integer), by the so-called {\it weighted counting},
$$
|A|_x \, := \,\sum_{a \in A} x^{f(a)} \quad,
$$
where $f:=A \rightarrow Z$ is the statistic in question. To go  from the weighted enumeration (a certain Laurent polynomial)
to straight enumeration, one sets $x=1$, i.e. $|A|_1 = |A|$.

Since this is {\it mathematics}, and not {\it accounting}, the usual scenario is not just {\bf one} specific set
$A$, but a sequence of sets $\{A_n\}_{n=0}^{\infty}$, and then the enumeration problem is to have an efficient
description of the numerical sequence $a_n:=|A_n|$, ready to be looked-up (or submitted) to the OEIS, and
its corresponding sequence of polynomials $P_n(x):=|A_n|_x$.

It often happens that the statistic $f$, defined on $A_n$, has a {\it scaled limiting distribution}.
In other words, if you draw a {\it histogram} of $f$ on $A_n$,, and do the
obvious {\it scaling}, they get closer and closer to a certain {\it continuous} curve,  as $n$ goes to infinity.

The scaling is as follows. Let $E_n(f)$ and $Var_n(f)$ the {\it expectation} and
{\it variance} of the statistic $f$ defined on $A_n$, and define the {\it scaled} random variable, for $a\in A_n$, by
$$
X_n (a):= \frac{f(a)- E_n(f)}{\sqrt{Var_n(f)}} \quad .
$$

If you draw the histograms of $X_n(a)$ for large $n$, they look practically the same, and converge to some {\it continuous limit}.

A famous example is {\it coin tossing}. If $A_n$ is $\{-1,1\}^n$, and $f(v)$ is the sum of $v$, then the limiting distribution
is the {\it bell shaped  curve} aka {\it standard normal distribution} aka {\it Gaussian distribution}.

As explained in [Z4], a purely finitistic approach to finding, and proving, a limiting scaled distribution, is via the
{\it method of moments}. Using {\it symbolic computation}, the computer can rigorously prove {\it exact} expressions
for as many moments as desired, and often (like in the above case, see [Z4]) find a recurrence for the
sequence of moments.  This enables one to identify the limits of the scaled moments with the moments of the continuous limit
(in the example of coin-tossing [and many other cases], $\frac{e^{-x^2/2}}{\sqrt{2\pi}}$, 
whose moments are famously $1 , 0 , 1\cdot 3, 0, 1 \cdot 3 \cdot 5, 0 , 1 \cdot 3 \cdot 5 \cdot 7, 0, \dots$) .
Whenever this is the case the discrete family of random variables is called {\it asymptotically normal}.
Whenever this is {\bf not} the case, it is interesting and surprising.

{\bf The Sum and Area Statistics on $a$-parking functions}

Let $\P(n,a)$ be the set of $a$-parking functions of length $n$.

A natural statistic is the sum
$$
Sum(p_1, \dots, p_n) := p_1 + p_2 + \dots + p_n =\sum_{i=1}^{n} p_i \quad .
$$
Another, even more natural (see the beautiful article [DH]) happens to be
$$
Area(p):= \frac{n(2a+n-1)}{2} -Sum(p) \quad .
$$

Let $P(n,a)(x)$ be the weighted analog of $p(n,a)$, according to {\it Sum}, i.e.
$$
P(n,a)(x) \, := \, \sum_{p \in  \P(n,a)} x^{Sum(p)} \quad .
$$

Analogously, let $Q(n,a)(x)$ be the weighted analog of $p(n,a)$, according to {\it Area}, i.e.
$$
Q(n,a)(x) \, := \, \sum_{p \in  \P(n,a)} x^{Area(p)} \quad .
$$

Clearly, one can easily go from one to the other
$$
Q(n,a)(x) \, = \, x^{(2a+n-1)n/2} \, P(n,a)(x^{-1}) \quad , \quad
P(n,a)(x) \, = \, x^{(2a+n-1)n/2} \, Q(n,a)(x^{-1}) \quad .
$$

How do we compute $P(n,a)(x)$?, (or equivalently, $Q(n,a)(x)$?). It is readily seen that the analog of
$(FundamentalRecurrence)$ for the weighted counting is
$$
P(n,a)(x) \, = \, x^n \, \sum_{k=0}^{n} \, {{n} \choose {k}} P(n-k,a+k-1)(x) \quad,
\eqno(FundamentalRecurrenceX)
$$
subject to the initial conditions $P(0,a)(x)=1$ and $P(n,0)(x)=0$.

So it is {\it almost} the same, the ``only" change is sticking $x^n$ in front of the sum on the right hand side.

Equivalently, 
$$
Q(n,a)(x) \, = \, \, \sum_{k=0}^{n} \, {{n} \choose {k}} x^{k(k+2a-3)/2} \,  Q(n-k,a+k-1)(x) \quad ,
\eqno(FundamentalRecurrenceAreaX)
$$
subject to the initial conditions $Q(0,a)(x)=1$ and $Q(n,0)(x)=0$.

Once again, in the sense of Wilf, this is already an {\it answer}, but because of the extra variable $x$, one
can not go as far as we did before for the naive, merely numeric, counting.

It is very unlikely that there is a ``closed form'' expression for $P(n,a)(x)$ (and hence $Q(n,a)(x)$), but for
{\it statistical purposes} it would be nice to get ``closed form'' expressions for

$\bullet$ the expectation,

$\bullet$ the variance,

$\bullet$  as many factorial moments as possible, from which the `raw' moments, and latter the {\it centralized} moments
and finally the {\it scaled moments} can be gotten. Then we can take the limits as $n$ goes to infinity, and
see if they match the moments of any of the known continuous distributions, and prove {\it rigorously}
that, at least for that many moments, the conjectured limiting distribution  matches.

In our case, the limiting distribution is the intriguing so-called {\it Airy distribution}, that Svante Janson prefers
to call ``area under Brownian excursion". This result was stated and proved in [DH], by using deep and sophisticated
{\it continuous} probability theory and continuous martingales. Here we will ``almost" prove this result, in the sense of
showing that the limits of the scaled moments of the area statistic on parking functions coincide
with the scaled moments of the Airy distribution up to the $30$-th moment, and we can go much further.

But we can do much more than continuous probabilists. We (or rather our computers, running Maple) can
find {\it exact} {\bf polynomial} expressions in $n$ and the expectation $E_1(n)$. 
We can do it for any desired number of moments, say $30$. Unlike continuous probability theorists,
our methods are entirely elementary, only using {\it high school algebra}.

We can also do the same thing for the more general $a$-parking functions. Now the expressions are
polynomials in $n$, $a$, and the expectation $E_1(n,a)$.

Finally, we believe that our approach, using the  recurrence $(FundamentalRecurrenceAreaX)$, can be used
to give a full proof (for all moments), by doing it asymptotically, and deriving a recurrence
for the leading terms of the asymptotics for the factorial moments that would coincide with the
well-known recurrence for the moments of the Airy distribution given, for example in  Eqs. (4) and (5) of
Svante Janson's article [J].  This is left as a challenge to our readers.

{\bf Finding the Expectation}

The expectation of the sum statistic, let's call it $E_{sum}(n,a)$ is given by (the prime denotes, as usual, differentiation w.r.t. $x$)
$$
E_{sum}(n,a) \, = \, \frac{P'(n,a)(1)}{P(n,a)(1)} \,  = \,  \frac{P'(n,a)(1)}{a(a+n)^{n-1}} \quad .
$$

Can we get a {\it closed-form expression} for $P'(n,a)(1)$, and hence for $E_{sum}(n,a)$?

Differentiating $(FundamentalRecurrenceX)$ with respect to $x$, using the product rule, we get
$$
P(n,a)'(x) \, =\, x^n \, \sum_{k=0}^{n} \, {{n} \choose {k}}  P(n-k,a+k-1)'(x) \, + \,
 \, n x^{n-1} \, \sum_{k=0}^{n} \, {{n} \choose {k}}  \, P(n-k,a+k-1)(x)  \quad .
$$
Plugging-in $x=1$ we get that $P(n,a)'(1)$, satisfies the recurrence 
$$
P(n,a)'(1) - \sum_{k=0}^{n} \, {{n} \choose {k}} P(n-k,a+k-1)'(1) =
n\, \sum_{k=0}^{n} \, {{n} \choose {k}} P(n-k,a+k-1)(1) = n\,p(n,a) \quad .
\eqno(FundamentalRecurrenceX1)
$$

Using this recurrence, we can, just as we did for $p(n,a)$ above, get  expressions, as polynomials in $a$,
for numeric $1 \leq n \leq 10$, say, and then conjecture that
$$
P'(n,a)(1)= \frac{1}{2} \, a\,n \, (a+n-1) \, (a+n)^{n-1}-
\frac{1}{2} \sum_{j=1}^{n} {{n} \choose {j}} \, j! \, a \, (a+n)^{n-j} \quad .
$$
To prove it, one plugs in the left side into $(FundamentalRecurrenceX1)$, changes the order of summation, and
simplifies. This is rather tedious, but since at the {\it end of the day}, these are equivalent to
{\it polynomial} identities in $n$ and $a$, checking it for sufficiently many special values of $n$ and $a$
would be a rigorous proof.

It follows that
$$
E_{sum}(n,a) \,= \,
\frac{n(a+n+1)}{2} - \frac{1}{2} \, \sum_{j=1}^{n} \frac{n!}{(n-j)! (a+n)^{j-1}} \quad .
$$
This formula first appears in [KY1].

Equivalently,

$$
E_{area}(n,a) \,= \,
\frac{n\,(a-2)}{2} \,+ \,\frac{1}{2} \, \sum_{j=1}^{n} \frac{n!}{(n-j)! (a+n)^{j-1}} \quad .
$$
In particular, for the primary object of interest, the case $a=1$, we get
$$
E_{area}(n,1) \,= \,
- \frac{n}{2} \,+ \,\frac{1}{2} \, \sum_{j=1}^{n} \frac{n!}{(n-j)! (n+1)^{j-1}} \quad .
$$
This rings a bell! It may written as
$$
E_{area}(n,1) \,= \,- \frac{n}{2} \,+ \, \frac{1}{2} W_{n+1} \quad,
$$
where $W_n$ is the {\bf iconic} quantity, 
$$
W_n \, = \, \frac{n!}{n^{n-1}} \sum_{k=0}^{n-2} \frac{n^k}{k!} \quad ,
$$
proved by Riordan and Sloane ([RS]) to be the expectation of another very important quantity, the sum of the heights
on rooted labeled trees on $n$ vertices. In addition to its considerable mathematical interest, this quantity, $W_n$,
has great {\it historical significance}, it was the {\it first sequence} , sequence $A435$
of the amazing On-Line Encyclopedia of Integer Sequences
(OEIS), now with almost  $300000$ sequences! See [EZ] for details, and far-reaching extensions, analogous to the present paper.

[{\eightrm The reason it is not sequence A1  is that initially the sequences were arranged in lexicographic order.}]

Another  fact, that will be of great use later in this paper, is that, as noted in [RS], Ramanujan and Watson proved that
 $W_n$ (and hence $W_{n+1}$) is asymptotic to 
$$
\frac{\sqrt{2\pi}}{2} \, n^{3/2} \quad .
$$

It is very possible that the formula $E_{area}(n,1) \,= \,- \frac{n}{2} \,+ \, \frac{1}{2} W_{n+1}$ may also be
deduced from the Riordan-Sloane result via one of the numerous known bijections between parking functions
and rooted labeled trees. More generally, the results below, for the special case $a=1$, might be deduced, from
those of [EZ], but we believe that the present {\it methodology} is interesting for its own sake, and besides
in our {\it current} approach (that uses recurrences rather than the Lagrange Inversion Formula),
it is much faster to compute higher moments, hence, going in the other direction, would produce many more moments for
the statistic on rooted labeled trees considered in [EZ], provided that there is indeed such a correspondence
that sends the area statistic on parking functions (suitably tweaked) to the Riordan-Sloane statistic on rooted labeled trees.

{\bf The Limiting Distribution}

Given a combinatorial family, one can easily get an idea of the limiting distribution
by taking a large enough $n$, say $n=100$, and generating a large enough number of random
objects, say $50000$, and drawing a {\it histogram}, see Figure 2 in Diaconis and
Hicks' insightful article [DH]. But, one does not have to resort to simulation.
While it is impractical to consider {\it all} $101^{99}$ parking functions of length $100$, the
generating function $Q(100,1)(x)$ contains the exact count for each
conceivable area from $0$ to ${{100} \choose {2}}$. See

{\tt http://sites.math.rutgers.edu/\~{}zeilberg/tokhniot/picsParking/Ha100.html} \quad ,

for the full histogram.

But an even more informative way to investigate the limiting distribution
 is to draw the histogram of the probability generating function
of the scaled distribution 
$$
X_n (p):= \frac{Area(p)- E_n}{\sqrt{Var_n}} \quad ,
$$
where $E_n$ and $Var_n$ are the expectation and variance respectively.

See

{\tt http://sites.math.rutgers.edu/\~{}zeilberg/tokhniot/picsParking/Da100.html } \quad ,

for $n=100$ and

{\tt http://sites.math.rutgers.edu/\~{}zeilberg/tokhniot/picsParking/Da120.html } \quad ,

for $n=120$. They look the same!

As proved in [DH] (using deep results in continuous probability due to David Aldous, Svante Janson,
and Chassaing and Marcket) the limiting distribution is the {\it Airy distribution}.
We will soon ``almost" prove it, but do much more by discovering exact expressions for the first $30$ moments,
not just their limiting asymptotics.

\vfill\eject

{\bf Truly Exact Expressions for the Factorial (and hence Centralized Moments)}

In [KY2] there is an ``exact'' expression for the general moment, that is not very useful for our purposes.
If one traces their proof, one can, conceivably, get explicit expressions for each specific moment,
but they did not bother to implement it, and the asymptotics is not immediate.

We discovered,  the following important  fact.

{\bf Fact.} Let $E_1(a,n):=E_{area}(a,n)$ be the expectation of the area statistic on $a$-parking functions
of length $n$,  given above, and let $E_k(n,a)$ be the $k$-th factorial moment
$$
E_k(n,a) \, := \, \frac{ Q^{(k)}(n,a)(1)}{a(a+n)^{n-1}} \quad ,
$$
then there exist polynomials $A_k(a,n)$ and $B_k(a,n)$ such that
$$
E_k(n,a) \, = A_k(a,n) \,+\,  B_k(a,n) \, E_1(a,n) \quad .
$$

The beauty of experimental mathematics is that these can be found by cranking out enough data, using
the sequence of probability generating functions $Q(n,a)(x)$, obtained by using  the recurrence,
$(FundamentalRecurrenceAreaX)$,  getting sufficiently many numerical data for the moments,
and using {\it undetermined coefficients}. These can be proved {\it a posteriori} by taking
these truly exact formulas and  verifying that the implied recurrences for the $k$-th factorial moment
(obtained from differentiating  $(FundamentalRecurrenceAreaX)$ $k$ times, using Leinitz's rule),
in terms of the previous ones. But this is {\bf not} necessary. Since, at the end of the day,
it all boils down to verifying {\bf polynomial identities}, so, once again, verifying them
for sufficiently many different values of $(n,a)$ constitutes a rigorous proof. To be fully rigorous,
one  needs to prove {\it a priori} bounds for the degrees in $n$ and $a$, but, in our humble opinion,
it is not that important, and could be left to the obtuse reader.

Our beloved computers, running the Maple package {\tt ParkingStatistics.txt}, available from the front of this article

{\tt http://sites.math.rutgers.edu/\~{}zeilberg/mamarim/mamarimhtml/par.html} \quad ,

produced the following, for the most interesting case of $a=1$, i.e. classical parking functions.

{\bf Theorem 1.} (equivalent to a result in [KY1]): The expectation of the area statistic on parking functions of length $n$ is
$$
E_1(n):= \,- \frac{n}{2} \,+ \, \frac{1}{2} \, \frac{(n+1)!}{(n+1)^{n}} \sum_{k=0}^{n-1} \frac{(n+1)^k}{k!} \quad ,
$$
and asymptotically it equals $\frac{\sqrt{2\pi}}{4} \cdot n^{3/2} +O(n)$.

{\bf Theorem 2.} The {\bf second} factorial moment of the area statistic on parking functions of length $n$ is
$$
 -\frac{7}{3} (n+1) \, E_1(n)+{\frac {5}{12}}\,{n}^{3}- \frac{1}{12} \,{n}^{2}- \frac{1}{3} \,n \quad ,
$$
and asymptotically it equals $\frac{5}{12} \cdot n^3 +O(n^{5/2})$.

{\bf Theorem 3.} The {\bf third} factorial moment of the area statistic on 
parking functions of length $n$ is
$$
-{\frac {175}{192}}\,{n}^{4}-{\frac {283}{192}}\,{n}^{3}+{\frac {199}{
192}}\,{n}^{2}+{\frac {259}{192}}\,n+ \left( {\frac {15}{32}}\,{n}^{3}+
{\frac {521}{96}}\,{n}^{2}+{\frac {1219}{96}}\,n+{\frac {743}{96}} 
 \right) \, E_1(n) \quad ,
$$
and asymptotically it equals $\frac{15}{128}\sqrt{2\pi} \cdot n^{9/2} +O(n^{4})$.

{\bf Theorem 4.} The {\bf fourth} factorial moment of the area statistic on 
parking functions of length $n$ is
$$
{\frac {221}{1008}}\,{n}^{6}+{\frac {63737}{30240}}\,{n}^{5}+{\frac {
101897}{15120}}\,{n}^{4}+{\frac {22217}{5040}}\,{n}^{3}-{\frac {1375}{
189}}\,{n}^{2}-{\frac {187463}{30240}}\,n
$$
$$
+ \left( -{\frac {35}{16}}\,{n
}^{4}-{\frac {449}{27}}\,{n}^{3}-{\frac {130243}{2520}}\,{n}^{2}-{
\frac {7409}{105}}\,n-{\frac {503803}{15120}} \right) \, E_1(n) \quad ,
$$
and asymptotically it equals $\frac{221}{1008}\cdot n^{6} +O(n^{11/2})$.

{\bf Theorem 5.} The {\bf fifth} factorial moment of the area statistic on 
parking functions of length $n$ is
$$
-{\frac {105845}{110592}}\,{n}^{7}-{\frac {2170159}{290304}}\,{n}^{6}-{
\frac {99955651}{3870720}}\,{n}^{5}-{\frac {30773609}{725760}}\,{n}^{4}
-{\frac {94846903}{11612160}}\,{n}^{3}+{\frac {24676991}{483840}}\,{n}^
{2}+{\frac {392763901}{11612160}}\,n
$$
$$
+ \left( {\frac {565}{2048}}\,{n}^{
6}+{\frac {1005}{128}}\,{n}^{5}+{\frac {9832585}{165888}}\,{n}^{4}+{
\frac {1111349}{5184}}\,{n}^{3}+{\frac {826358527}{1935360}}\,{n}^{2}+{
\frac {159943787}{362880}}\,n+{\frac {1024580441}{5806080}} \right) \, E_1(n) \quad ,
$$
and asymptotically it equals $\frac{565}{8192    }\sqrt{2\pi} \cdot n^{15/2} +O(n^{7})$.

{\bf Theorem 6.} The {\bf sixth} factorial moment of the area statistic 
parking functions of length $n$ is
$$
{\frac {82825}{576576}}\,{n}^{9}+{\frac {373340075}{110702592}}\,{n}^{8
}+{\frac {9401544029}{332107776}}\,{n}^{7}+{\frac {14473244813}{
127733760}}\,{n}^{6}+{\frac {414139396709}{1660538880}}\,{n}^{5}
$$
$$
+{
\frac {88215445651}{332107776}}\,{n}^{4}-{\frac {18783816473}{332107776
}}\,{n}^{3}-{\frac {643359542029}{1660538880}}\,{n}^{2}-{\frac {
358936540409}{1660538880}}\,n
$$
$$
+ ( -{\frac {3955}{2048}}\,{n}^{7}-{
\frac {186349}{6144}}\,{n}^{6}-{\frac {259283273}{1161216}}\,{n}^{5}-{
\frac {119912501}{129024}}\,{n}^{4}-{\frac {149860633081}{63866880}}\,{
n}^{3}
$$
$$
-{\frac {601794266581}{166053888}}\,{n}^{2}-{\frac {864000570107}
{276756480}}\,n-{\frac {921390308389}{830269440}} ) \, E_1(n) \quad ,
$$
and asymptotically it equals $\frac{ 82825}{576576}\cdot n^{9} +O(n^{17/2})$.

For Theorems 7-30, see the output file

{\tt http://sites.math.rutgers.edu/\~{}zeilberg/tokhniot/oParkingStatistics7.txt} .

Let $\{e_k\}_{k=1}^{\infty}$ be the sequence of moments of the Airy distribution, defined by the recurrence given
in Equations $(4)$ and $(5)$ in Svante Janson's interesting survey paper [J]. Our computers, using our
Maple package, proved that 
$$
E_k(n) \, = \, e_k n^{\frac{3k}{2}} + O( n^{\frac{3k-1}{2}}) \quad,
$$
for $1 \leq k \leq 30$. It follows that the {\it limiting distribution} of the area statistic is (most probably) the Airy distribution,
since the first $30$ moments match.
Of course, this was already known to continuous probability theorists, and we only proved it for the first $30$
moments, {\bf but}:

$\bullet$ Our methods are purely {\it elementary} and {\it finitistic}

$\bullet$ We can easily go much farther, i.e. prove it for more moments

$\bullet$ We believe that our approach, using recurrences, can be used to derive a recurrence for the {\it leading}
asymptotics of the factorial moments, $E_k(n)$, that would turn out to be the same as the above mentioned
recurrence (Eqs. (4) and (5) in [J]). We leave this as a challenge to the reader.

{\bf Exact expressions for the first $10$ moments of the Area statistic for general $a$-parking}

To see expressions in $a$, $n$, and $E_1(n,a)$, for the first $10$ moments of $a$-parking, see

{\tt http://sites.math.rutgers.edu/\~{}zeilberg/tokhniot/oParkingStatistics8.txt} \quad .

{\bf Acknowledgment}: Many thanks are due to Valentin F\'eray and Svante Janson for insightful information and
useful references. Also thanks to Benji Weiss for comments on a previous version.

{\bf References}

[DH] Persi Diaconis and Angela Hicks, {\it Probabilizing Parking Functions}, Adv. in Appl. Math. {\bf89} (2017), 125-155. \hfill\break
{\tt https://arxiv.org/abs/1611.09821} \quad .

[EZ] Shalosh B. Ekhad and Doron Zeilberger,
{\it Going Back to Neil Sloane's FIRST LOVE (OEIS Sequence A435): On the Total Heights in Rooted Labeled Trees},
The Personal Journal of Shalosh B. Ekhad and Doron Zeilberger,  July 19, 2016. \hfill\break
{\tt http://sites.math.rutgers.edu/\~{}zeilberg/mamarim/mamarimhtml/a435.html} (accessed May 24, 2018) \quad .

[FR] Dominique Foata and John Riordan, {\it Mapping of acyclic and parking functions}, Aequationes Mathematicae {\bf 10} (1974),
490-515.

[J] Svante Janson, {\it Brownian excursion area, Wright’s constants in graph enumeration, and other Brownian areas}, Probab. Surveys,
{\bf 4}(2007), 80-145. \hfill\break
{\tt https://projecteuclid.org/euclid.ps/1178804352} (accessed May 24, 2018) \quad .

[KW] Alan G. Konheim and Benjamin Weiss, {\it An occupancy discipline and applications},
SIAM J. Applied Math. {\bf 14} (1966), 1266-1274. [Available from JSTOR.] 

[KY1] J.P. Kung and C. Yan, {\it Expected sums of general parking functions}, Annals of Combinatorics {\bf 7} (2003), 481-493.

[KY2] J.P. Kung and C. Yan, {\it Exact formulas for the moments of sums of classical parking functions}, Advances in Applied Mathematics {\bf 31}
(2003), 215-241.

[RS] John Riordan and Neil J. A. Sloane, 
{\it The enumeration of rooted trees by total height}, J. Australian Math. Soc. {\bf 10} (1969), 278-282.  \hfill\break
{\tt http://neilsloane.com/doc/riordan-enum-trees-by-height.pdf} \quad  (accessed May 24, 2018) .

[Sc] Marcel-Paul Sch\"utzenberger,  {\it  On an enumeration problem},
J. Combinatorial Theory  {\bf 4} (1968), 219-221.

[St1] Richard Stanley, {\it Parking functions}, \hfill\break
{\tt www-math.mit.edu/\~{}rstan/transparencies/parking.pdf} (accessed May 24, 2018) \quad .

[St2] Richard Stanley, {\it A survey of parking functions}, \hfill\break
{\tt www-math.mit.edu/\~{}rstan/transparencies/parking3.pdf} (accessed May 24, 2018) \quad .

[St3] Richard Stanley, {\it ``Enumerative Combinatorics, Volume 2''}, Cambridge University Press, 1999.

[W] Herbert S. Wilf, {\it  What is an Answer?}, The American Mathematical Monthly {\bf 89} (1982), 289-292.

[Z1] Doron Zeilberger,  {\it Symbolic Moment Calculus I.: Foundations and Permutation Pattern Statistics},
Annals of Combinatorics {\bf 8} (2004), 369-378.
\hfill\break
{\tt http://sites.math.rutgers.edu/\~{}zeilberg/mamarim/mamarimhtml/smcI.html} (accessed May 24, 2018) \quad .

[Z2] Doron Zeilberger, {\it Symbolic Moment Calculus II.: Why is Ramsey Theory Sooooo Eeeenormously Hard?},
INTEGERS {\bf 7(2)}(2007), A34.
\hfill\break
{\tt http://sites.math.rutgers.edu/\~{}zeilberg/mamarim/mamarimhtml/smcII.html} (accessed May 24, 2018) \quad .

[Z3] Doron Zeilberger,  {\it The $n^{n-2}$-th proof for the number of labeled trees}, The
Personal Journal of Shalosh B. Ekhad and Doron Zeilberger, undated (c. 1998), \hfill\break
{\tt http://sites.math.rutgers.edu/\~{}zeilberg/mamarim/mamarimhtml/labtree.html} \quad (accessed May 24, 2018).

[Z4] Doron Zeilberger,  {\it  The Automatic Central Limit Theorems Generator (and Much More!)}, 
``Advances in Combinatorial Mathematics: Proceedings of the Waterloo Workshop in Computer Algebra 2008 
in honor of Georgy P. Egorychev", chapter 8, pp. 165-174, (I.Kotsireas, E.Zima, eds. Springer Verlag, 2009.) \hfill\break
{\tt http://sites.math.rutgers.edu/\~{}zeilberg/mamarim/mamarimhtml/georgy.html} (accessed May 24, 2018) \quad .

\bigskip
\hrule
\bigskip
Yukun Yao, Department of Mathematics, Rutgers University (New Brunswick), Hill Center-Busch Campus, 110 Frelinghuysen
Rd., Piscataway, NJ 08854-8019, USA. \hfill\break
Email: {\tt yao at math dot rutgers dot edu}   \quad .
\bigskip
Doron Zeilberger, Department of Mathematics, Rutgers University (New Brunswick), Hill Center-Busch Campus, 110 Frelinghuysen
Rd., Piscataway, NJ 08854-8019, USA. \hfill\break
Email: {\tt DoronZeil at gmail  dot com}   \quad .
\bigskip
\hrule
\bigskip
First Written: June 5, 2018 ; This version: June 6, 2018.

\end